\documentstyle[amsfonts,12pt]{article}
\input epsf

\def\underscore#1{\underline{\vphantom{q}#1}}
\def\mapsisoto{
   \makebox[0pt][l]{\raisebox{1.75pt}{\hspace*{3.2pt}$\cong$}}\longrightarrow}
\def\rom#1{{\normalshape #1}}
\def\rank{\mathop{\mathrm{rank}}}
\newtheorem{thm}{Theorem}
\newtheorem{cor}{Corollary}

\newenvironment{prf}{\begin{trivlist}\item[]{\bf Proof} }%
{\hfill $\Box$ \end{trivlist}}
\begin{document}
\begin{center}\Large\bf Drawing with Complex Numbers\\[15pt]
\large\mediumseries Michael Eastwood\footnote{Supported by the Australian
Research Council.} and Roger Penrose
\end{center}

It is not commonly realized that the algebra of complex numbers can be used in
an elegant way to represent the images of ordinary $3$-dimensional figures,
orthographically projected to the plane. We describe these ideas here, both
using simple geometry and setting them in a broader context.

Consider orthogonal projection in Euclidean $n$-space onto an 
$m$-di\-men\-sion\-al
subspace. We may as well choose co\"ordinates so that this is the standard
projection $P:{\Bbb R}^n\to{\Bbb R}^m$ onto the first $m$ variables. Fix a
non-degenerate simplex $\Sigma$ in~${\Bbb R}^n$. Two such simplices are said to
be \underscore{similar} if one can be obtained from the other by a Euclidean
motion together with an overall scaling. This article answers the
following question. Given $n+1$ points in ${\Bbb R}^m$, when can these points 
be obtained as the images under $P$ of the vertices of a simplex similar
to~$\Sigma$?

When $n=3$ and $m=2$, then $P$ is the standard \underscore{orthographic}
projection (as often used in engineering drawing) and we are concerned with how
to draw a given tetrahedron. We shall show, for example, that four points
$\alpha,\beta,\gamma,\delta$ in the plane are the orthographic projections of
the vertices of a \underscore{regular} tetrahedron if and only if
\begin{equation}\label{tetrahedron}
(\alpha+\beta+\gamma+\delta)^2=4(\alpha^2+\beta^2+\gamma^2+\delta^2)
\end{equation}
where $\alpha,\beta,\gamma,\delta$ are regarded as \underscore{complex}
numbers! Similarly, suppose a cube is orthographically projected and normalised
so that a particular vertex is mapped to the origin. If $\alpha,\beta,\gamma$
are the images of the three neighbouring vertices, then
\begin{equation}\label{cube}\alpha^2+\beta^2+\gamma^2=0,\end{equation}
again as a \underscore{complex} equation. Conversely, if this equation is
satisfied, then one can find a cube whose orthographic image is given in this
way. Since parallel lines are seen as parallel in the drawing,
equation~(\ref{cube}) allows one to draw the general cube:
\begin{center}
\begin{picture}(100,100)(-20,0)
\put(-150,-50){\epsfysize=700pt \epsffile{gauss.eps}}
\put(153,97){$\gamma$}
\put(37,87){$\beta$}
\put(121,5){$\alpha$}
\put(116,74){$0$}
\put(-150,50){\begin{tabular}l In this example,  $\alpha=2-26i$\\
                      \phantom{In this example,} $\beta=-23+2i$\\
                      \phantom{In this example,} $\gamma=14+7i$\end{tabular}}
\end{picture}
\end{center}

The result for a cube is known as \underscore{Gauss'} \underscore{fundamental}
\underscore{theorem} \underscore{of} \linebreak \underscore{axonometry}---see
\cite[p.~309]{gauss} where it is stated without proof. In engineering drawing,
one usually fixes three \underscore{principal} axes in Euclidean three-space
and then an orthographic projection onto a plane transverse to these axes is
known as an \underscore{axonometric} projection (see, for example,
\cite[Chapter~17]{hoelscherandspringer}). Gauss' theorem may be regarded as
determining the degree of foreshortening along the principal axes for a general
axonometric projection. The projection corresponding to taking
$\alpha,\beta,\gamma$ to be the three cube roots of unity is called
\underscore{isometric} projection because the foreshortening is the same for
the three principal axes. In an axonometric drawing, it is conventional to take
the image axes at mutually obtuse angles:
\begin{center}
\begin{picture}(100,110)(50,5)
\put(-150,-50){\epsfysize=700pt \epsffile{weisbach.eps}}
\put(102,102){$A$}
\put(133.5,58.5){$B$}
\put(87.5,51.5){$C$}
\end{picture}
\end{center}
If $|\alpha|=a$, $|\beta|=b$, $|\gamma|=c$, then equation~(\ref{cube}) is
equivalent to the sine rule for the triangle with sides $\alpha^2$, $\beta^2$,
$\gamma^2$, namely
$$\frac{a^2}{\sin 2A}=\frac{b^2}{\sin 2B}=\frac{c^2}{\sin 2C}.$$
In this form, the fundamental theorem of axonometry is due to Weisbach and was
published in T\"ubingen in 1844 in the Polytechnische Mitteilungen of Volz and
Karmasch. Equivalent statements can be found in modern engineering drawing
texts (e.g.~\cite[p.~44]{rothandvanhaeringen}).

Equation (\ref{cube}) may be used to give a ruler and compass construction of
the general orthographic image of a cube. If we suppose that the image of a
vertex and two of its neighbours are already specified, then (\ref{cube})
determines (up to a two-fold ambiguity) the image of the third neighbour. The
construction is straightforward except perhaps for the construction of a 
complex
square root for which we advocate the following as quite efficient:
\begin{center}
\begin{picture}(100,110)(50,5)
\put(-150,-50){\epsfysize=700pt \epsffile{sqrt.eps}}
\put(132,104){$z$}
\put(153,73){$\sqrt z$}
\put(104,51){$0$}
\put(153,44){$1$}
\put(42,44){$\zeta$}
\put(109,45){$\bullet$}
\put(146,45){$\bullet$}
\put(146,73){$\bullet$}
\put(125,101){$\bullet$}
\put(50,45){$\bullet$}
\end{picture}
\end{center}
Firstly $\zeta$ is constructed by marking the real axis at a distance~$\|z\|$
from the origin. Then, a circle is constructed passing through the three points
$\zeta$, $1$, and~$z$. Finally, the angle between $1$ and $z$ is bisected and
$\sqrt z$ appears where this bisector meets the circle.

In engineering drawing, it is more usual that the images of the three principal
axes are prescribed or chosen by the designer and one needs to determine the
relative degree of foreshortening along these axes. There is a ruler and
compass construction given by T.~Schmid in 1922 (see, for example,
\cite[\S17.17--17.19]{hoelscherandspringer}):
\begin{center}
\begin{picture}(100,80)(50,5)
\put(-88,-177){\epsfysize=1200pt \epsffile{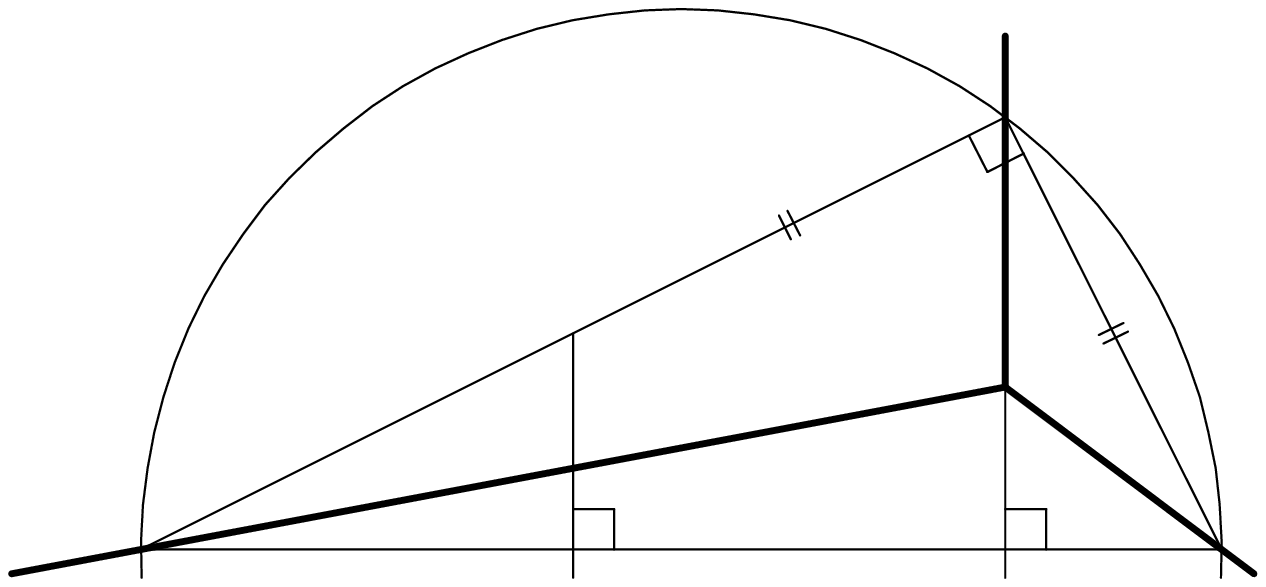}}
\put(177,14){$\alpha$}
\put(91,30){$\beta$}
\put(20,17){$P$}
\put(149,70){$Q$}
\put(75,43){$R$}
\put(134,37){$0$}
\put(170.3,11){$\bullet$}
\put(86,22){$\bullet$}
\put(29.5,11){$\bullet$}
\put(142,67){$\bullet$}
\put(86,39){$\bullet$}
\put(142,32.3){$\bullet$}
\end{picture}
\end{center}
In this diagram, the three principal axes and $\alpha$ are given. By drawing a
perpendicular from $\alpha$ to one of the of the principal axes and marking its
intersection with the remaining principal axis, we obtain~$P$. The point $Q$ is
obtained by drawing a semi-circle as illustrated. The point $R$ is on the
resulting line and equidistant with $\alpha$ from~$Q$. Finally, $\beta$ is
obtained by dropping a perpendicular as shown. It is easy to see that this
construction has the desired effect---in Euclidean three-space rotate the
right-angled triangle with hypotenuse $P\alpha$ about this hypotenuse until the
point $Q$ lies directly above $0$ in which case $R$ will lie directly
above~$\beta$ and the third vertex will lie somewhere over the line through $0$
and~$Q$. One may verify the appropriate part of Weisbach's condition
\begin{equation}\label{partialweisbach}
\frac{a^2}{\sin 2A}=\frac{b^2}{\sin 2B}\end{equation}
by the following calculation. Without loss of generality we may represent all
these points by complex numbers normalised so that~$Q=1$. Then it is
straightforward to check that
$$R\!=\!1+i-i\alpha,\;
P\!=\!\frac{\alpha(\alpha+\overline\alpha)+2(1-\alpha-\overline\alpha)}
       {\alpha-\overline\alpha},\;
\beta\!=\!\frac{\alpha(\alpha+\overline\alpha)+2(1-\alpha-\overline\alpha)}
           {2-\alpha-\overline\alpha}i$$
and therefore that
$$\alpha^2+\beta^2 =
  4\frac{(\alpha-1)\overline{(\alpha-1)}(\alpha+\overline\alpha-1)}
        {(\alpha+\overline\alpha-2)^2}.$$
That $\alpha^2+\beta^2$ is real is equivalent to~(\ref{partialweisbach}).

To prove Gauss' theorem more directly consider three vectors in ${\Bbb R}^3$ 
as the columns of a $3\times 3$ matrix. This matrix is orthogonal
if and only if the three vectors are orthonormal. It is
equivalent to demand that the three rows be orthonormal.
However, any two orthonormal vectors in ${\Bbb R}^3$ may be extended
to an orthonormal basis. Thus, the condition that three vectors
$$\left\lgroup\begin{array}c x_1\\y_1\end{array}\right\rgroup\qquad
\left\lgroup\begin{array}c x_2\\y_2\end{array}\right\rgroup\qquad
\left\lgroup\begin{array}c x_3\\y_3\end{array}\right\rgroup$$
in ${\Bbb R}^2$ be the images under $P:{\Bbb R}^3\to{\Bbb R}^2$ of an
orthonormal basis of ${\Bbb R}^3$, is that
$$\left\lgroup\begin{array}{ccc} x_1&x_2&x_3\end{array}\right\rgroup
\quad\mbox{ and }\quad
\left\lgroup\begin{array}{ccc} y_1&y_2&y_3\end{array}\right\rgroup$$
be orthonormal in~${\Bbb R}^3$. Dropping the overall scale, we obtain
$$x_1{}^2+x_2{}^2+x_3{}^2=y_1{}^2+y_2{}^2+y_3{}^2\quad\mbox{ and }\quad
x_1y_1+x_2y_2+x_3y_3=0.$$
Writing, $\alpha=x_1+iy_1$, $\beta=x_2+y_2$, $\gamma=x_3+y_3$, these two
equations are the real and imaginary parts of~(\ref{cube}). To deduce the case
of a regular tetrahedron as described by equation~(\ref{tetrahedron}) from the
case of a cube as described by equation~(\ref{cube}), it suffices to note that
equation~(\ref{tetrahedron}) is translation invariant and that a regular
tetrahedron may be inscribed in a cube. Thus, we may take
$\delta=\alpha+\beta+\gamma$ and observe that (\ref{tetrahedron}) and
(\ref{cube}) are then equivalent.

It is easy to see that the possible images of a particular tetrahedron $\Sigma$
in ${\Bbb R}^3$ under an arbitrary Euclidean motion followed by the projection
$P$ form a $5$-dimensional space---the group of Euclidean motions is
$6$-dimensional but translation orthogonal to the plane leaves the image
unaltered. It therefore has codimension $3$ in the $8$-dimensional space of all
tetrahedral images ($2$~degrees of freedom for each vertex). Allowing similar
tetrahedra rather than congruent reduces the codimension to~$2$. Therefore, two
real equations are to be expected. Always, these two real equations combine as
a single \underscore{complex} equation such as (\ref{tetrahedron})
or~(\ref{cube}). At first sight, this is perhaps surprising and even more so
when the same phenomenon occurs for $P:{\Bbb R}^n\to{\Bbb R}^2$ for 
arbitrary~$n$.

For $n=3$, there is a proof of Gauss' theorem which brings in complex numbers
at the outset. Consider the space $H$ of Hermitian $2\times 2$ matrices with
zero trace, i.e. matrices of the form
$$X=\left\lgroup\begin{array}{cc}w&u+iv\\ u-iv&-w\end{array}\right\rgroup\quad
\mbox{for }
\left\lgroup\begin{array}c u\\ v\\ w\end{array}\right\rgroup\in{\Bbb R}^3.$$
We may identify $H$ with ${\Bbb R}^3$ and, in so doing, $-\det X$ becomes the
square of the Euclidean length. The group $G$ of invertible $2\times 2$ complex
matrices of the form
$$\Lambda=\left\lgroup\matrix{a & -b\cr
\overline b & \overline a}\right\rgroup$$
acts linearly on $H$ by $X\mapsto\Lambda X\overline\Lambda^t$. Moreover,
$$\det(\Lambda X\overline\Lambda^t)=(|a|^2+|b|^2)^2\det X$$
so $G$ acts by similarities. It is easy to check that all similarities may be
obtained in this way. (This trick is essentially as used in Hamilton's theory
of quaternions and is well known to physicists---in modern parlance it is
equivalent to the isomorphism of Lie groups
${\mathrm{Spin}}(3)\cong{\mathrm{SU}}(2)$.) Therefore, an arbitrary 
orthographicimage of a cube may be obtained by acting with $\Lambda$ on the 
standard basis
$$\left\lgroup\begin{array}{cc}0&1\\ 1&0\end{array}\right\rgroup,
\quad\left\lgroup\begin{array}{cc}0&i\\ -i&0\end{array}\right\rgroup,
\quad\left\lgroup\begin{array}{cc}1&0\\ 0&-1\end{array}\right\rgroup$$
and then picking out the top right hand entries. We obtain
$$\begin{array}{rclclcl}
\Lambda\left\lgroup\begin{array}{cc}0&1\\ 1&0\end{array}\right\rgroup
\overline\Lambda^t & =
&\left\lgroup\begin{array}{cc}\ast&a^2-b^2\\ \ast&\ast\end{array}\right\rgroup
&\longmapsto& a^2-b^2&=&\alpha\\[20pt]
\Lambda\left\lgroup\begin{array}{cc}0&i\\ -i&0\end{array}\right\rgroup
\overline\Lambda^t & =
&\left\lgroup\begin{array}{cc}\ast&i(a^2+b^2)\\
                 \ast&\ast\end{array}\right\rgroup
&\longmapsto& i(a^2+b^2)&=&\beta\\[20pt]
\Lambda\left\lgroup\begin{array}{cc}1&0\\ 0&-1\end{array}\right\rgroup
\overline\Lambda^t & =
&\left\lgroup\begin{array}{cc}\ast&2ab\\ \ast&\ast\end{array}\right\rgroup
&\longmapsto& 2ab&=&\gamma\end{array}$$
and therefore $\alpha^2+\beta^2+\gamma^2=0$, as required. Conversely,
this is exactly the condition that $\alpha,\beta,\gamma$ may be written in this
form. (Compare the half angle formulae---if $s^2+c^2=1$, then $s=2t/(1+t^2)$
and $c=(1-t^2)/(1+t^2)$ for some~$t$.) That Gauss~\cite[p.~309]{gauss} makes
the same observation concerning the form of $\alpha,\beta,\gamma$ suggests that
perhaps he also had this reasoning in mind.

The proof of Gauss' theorem using orthogonal matrices clearly extends to
$P:{\Bbb R}^n\to{\Bbb R}^2={\Bbb C}$ for arbitrary~$n$. To state it, the 
following 
terminology concerning the standard projection $P:{\Bbb R}^n\to{\Bbb R}^m$ is
useful. We shall say that $v_1,v_2,\ldots,v_n\in{\Bbb R}^m$ are
\underscore{normalised} \underscore{eutactic} if and only if there is an
orthonormal basis $u_1,u_2,\ldots,u_n$ of ${\Bbb R}^n$ with $v_j=Pu_j$
for~$j=1,2,\ldots,n$. We shall say that $v_1,v_2,\ldots,v_n\in{\Bbb R}^m$ are
\underscore{eutactic} if and only if $\mu v_1,\mu v_2,\ldots,\mu v_n$ are
normalised eutactic for some~$\mu\not=0$.
\begin{thm} The points $z_1,z_2,\ldots,z_n\in{\Bbb C}={\Bbb R}^2$ are 
eutactic if and only if
$$z_1{}^2+z_2{}^2+\cdots+z_n{}^2=0$$
and not all $z_j$ are zero.\end{thm}
There is a proof for $n=4$ based on the isomorphism
$${\mathrm{Spin}}(4)\cong{\mathrm{SU}}(2)\times{\mathrm{SU}}(2)$$
and, indeed, this is how we came across the theorem in the first place.
However, a more direct route to complex numbers and one which applies in all
dimensions is based on the observation that ${\mathrm{Gr}}_2^+({\Bbb R}^2)$, 
the Grassmannian of oriented two-planes in ${\Bbb R}^n$, is naturally a
\underscore{complex} manifold. When $n=3$, this Grassmannian is just the
two-sphere and has a complex structure as the Riemann sphere. In general,
consider the mapping
$${\Bbb{CP}}_{n-1}\setminus{\Bbb{RP}}_{n-1}\stackrel{\pi}{\longrightarrow}
                                              {\mathrm{Gr}}_2^+({\Bbb R}^n)$$
induced by~${\Bbb C}^n\ni z\mapsto iz\wedge\overline z$. In other words, a
complex vector $z=x+iy\in{\Bbb C}^n$ is mapped to the two-dimensional oriented
subspace of ${\Bbb R}^n$ spanned by $x$ and $y$, the real and imaginary parts
of~$z$. Let $\langle\phantom{z},\phantom{z}\rangle$ denote the standard inner
product on ${\Bbb R}^n$ extended to ${\Bbb C}^n$ as a complex bilinear form. 
Then, $\langle z,z\rangle=0$ imposes two real equations
$$\|x\|^2=\|y\|^2\quad\mbox{and}\quad\langle x,y\rangle=0$$
on the real and imaginary parts. In other words, $x,y$ is proportional to an
orthonormal basis for~${\mathrm{span}}\{x,y\}$. Hence, if $z$ and $w$ satisfy
\mbox{$\langle z,z\rangle=0=\langle w,w\rangle$} and define the same oriented
two-plane, then $w=\lambda z$ for some~$\lambda\in{\Bbb C}\setminus\{0\}$. The
non-singular complex quadric
$$K=\{[z]\in{\Bbb{CP}}_{n-1}\mbox{ s.t. }\langle z,z\rangle=0\}$$
avoids ${\Bbb{RP}}_{n-1}\subset{\Bbb{CP}}_{n-1}$ and we have shown that 
$\pi|_K$ is injective. It is clearly surjective. The isomorphism
$$\pi:K\mapsisoto{\mathrm{Gr}}_2^+({\Bbb R}^n)$$
respects the natural action of ${\mathrm{SO}}(n)$ on $K$
and~${\mathrm{Gr}}_2^+({\Bbb R}^n)$. The generalised Gauss theorem follows
immediately since, rather than asking about the image of a general orthonormal
basis under the standard projection $P:{\Bbb R}^n\to{\Bbb R}^2$, we may,
equivalently, ask about the image of the standard basis $e_1,e_2,\ldots,e_n$
under a general orthogonal projection onto an oriented
two-plane~$\Pi\subset{\Bbb R}^n$. Any such $\Pi$ is naturally complex, the 
action of $i$ being given by rotation by $90^\circ$ in the positive sense. 
If $\Pi$ is represented by $[z_1,z_2,\ldots,z_n]\in K$ as above and we use 
$x,y\in\Pi$ to identify $\Pi$ with~${\Bbb C}$, then $e_j\mapsto z_j$ and
$$z_1{}^2+z_2{}^2+\cdots+z_n{}^2=\langle z,z\rangle=0,$$
as required. Conversely, a solution of this complex equation determines an
appropriate plane~$\Pi$.

For the case of a general tetrahedron or simplex and for general $m$ and $n$,
it is more convenient to start with Hadwiger's theorem \cite{hadwiger} or
\cite[page 251]{coxeter} as follows. The proof is obtained by extending our
orthogonal matrix proof of Gauss' theorem.
\begin{thm}[Hadwiger] Assemble $v_1,v_2,\ldots,,v_n\in{\Bbb R}^m$ as the 
columns of an $m\times n$ matrix~$V$. These vectors are normalised eutactic 
if and only if~$VV^t=1$ \rom{(}the $m\times m$ identity matrix\rom{)}.
\end{thm}
\begin{prf} If $v_1,v_2,\ldots,v_n$ are normalised eutactic, then assembling a
corresponding orthonormal basis of ${\Bbb R}^n$ as the columns of an 
$n\times n$ matrix, we have $V=PU$ and $U^tU=1$ (the $n\times n$ identity 
matrix). Therefore, $UU^t=1$ and
$$VV^t=PUU^tP^t=PP^t=1,$$
as required. Conversely, if $VV^t=1$, then the columns of $V^t$ may be
completed to an orthonormal basis of ${\Bbb R}^n$, i.e.\ $V^t=U^tP^t$ for
$UU^t=1$. Now, $U^tU=1$ and $V=PU$, as required. \end{prf}

The case of a general simplex is obtained essentially by a change of basis as
follows. Suppose $a_1,a_2,\ldots,a_n,a_{n+1}$ are the vertices of a
non-degenerate simplex $\Sigma$ in~${\Bbb R}^n$ whose centre of mass is at the
origin. In other words, the $n\times(n+1)$ matrix $A$ has rank $n$ and $Ae=0$
where $e$ is the column vector all of whose $n+1$ entries are~$1$. Form the
$(n+1)\times(n+1)$ symmetric matrix
$$Q=A^t(AA^t)^{-2}A,$$
noting that $\rank A=n$ implies the \underscore{moment} \underscore{matrix}
$AA^t$ is invertible.
\begin{thm} Given $b_1,b_2,\ldots,b_n,b_{n+1}\in{\Bbb R}^m$ assembled as the
columns of an $m\times(n+1)$ matrix $B$, these vectors are the images under
orthogonal projection $P:{\Bbb R}^n\to{\Bbb R}^m$ of the vertices of a simplex
congruent to $\Sigma$ if and only if
\begin{equation}\label{simplex} BQ^tB=1.\end{equation}\end{thm}
\begin{prf} The vertices of a simplex congruent to $\Sigma$ are the columns of
a matrix $UA+ae^t$ for some orthogonal matrix $U$ and translation vector
$a\in{\Bbb R}^n$. Also, note that $Qe=0$. Thus, if $B=P(UA+ae^t)$, then
$$\begin{array}{rcl}
BQB^t&=&PUAQA^tU^tP^t\\
     &=&PUAA^t(AA^t)^{-2}AA^tU^tP^t\\
     &=&PUU^tP^t\,=\,PP^t\,=\,1,\end{array}$$
as required. Conversely, $Qe=0$ implies that (\ref{simplex}) is translation
invariant. So, without loss of generality, we may suppose that
$b_1+b_2+\cdots+b_n+b_{n+1}=0$, that is to say, $Be=0$. Writing out
(\ref{simplex}) in full gives
$$BA^t(AA^t)^{-1}(BA^t(AA^t)^{-1})^t=1$$
so, by Hadwiger's theorem, there is an orthogonal matrix $U$ so that
$$BA^t(AA^t)^{-1}=PU.$$
Thus,
$$BA^t(AA^t)^{-1}A=PUA\quad\mbox{ and }\quad Be=0.$$
Certainly, $B=PUA$ is a solution of these equations but it is the only solution
since $A^t(AA^t)^{-1}A$ has rank $n$ and $e$ is not in the range of this linear
transformation. \end{prf}
\begin{cor}[case $\bold{m=2}$] Points
$z_1,z_2,\ldots,z_n,z_{n+1}\in{\Bbb C}$ are the images under orthogonal
projection of the vertices of a simplex similar to $\Sigma$ if and only if
$$z^tQz=0$$
where $z$ is the column vector with components
$z_1,z_2,\ldots,z_n,z_{n+1}$.\end{cor}

It is, of course, possible to compute $Q$ explicitly for any given example. If
the simplex $\Sigma$ has some degree of symmetry, however, we can often
circumvent such computation. Consider, for example, the case of a
\underscore{regular} simplex. From the corollary above, we know that the image
of such a simplex in the plane is characterised by a complex homogeneous
quadratic polynomial. The symmetries of the regular simplex ensure that this
polynomial must be invariant under ${\cal S}_{n+1}$, the symmetric group on 
$n+1$ letters. Hence, it must be expressible in terms of the elementary 
symmetric polynomials. Equivalently, it must be a linear combination of
$$(z_1+z_2+\cdots+z_n+z_{n+1})^2\quad\mbox{ and }\quad
   z_1{}^2+z_2{}^2+\cdots+z_n{}^2+z_{n+1}{}^2.$$
Up to scale, there is only one such combination which is translation invariant,
namely
\begin{equation}\label{regularsimplex}(z_1+z_2+\cdots+z_n+z_{n+1})^2-
  (n+1)(z_1{}^2+z_2{}^2+\cdots+z_n{}^2+z_{n+1}{}^2).\end{equation}
It follows that the vanishing of this polynomial is an equation which
characterises the possible images of a regular simplex under orthogonal
projection into the plane. The special case $n=2$ characterises the equilateral
triangles in the plane \cite[Problem 15 on page 79]{barnardandchild}.

Equation (\ref{cube}) characterising the orthographic images of a cube, may be
deduced by similar symmetry considerations. If a particular vertex is mapped to
the origin and its neighbours are mapped to $\alpha,\beta,\gamma$ then, since
each of these neighbouring vertices is on an equal footing, the polynomial in
question must be a linear combination of $(\alpha+\beta+\gamma)^2$ and
$\alpha^2+\beta^2+\gamma^2$. To find out which linear combination we need only
consider a particular projection, for example:
\begin{center}
\begin{picture}(70,70)
\put(10,60){\line(1,0){50}}
\put(60,10){\line(0,1){50}}
\put(10,10){\makebox(0,0){$\bullet$}}
\put(60,10){\makebox(0,0){$\bullet$}}
\put(10,60){\makebox(0,0){$\bullet$}}
\put(-7,67){\makebox(0,0){$\gamma=i$}}
\put(77,3){\makebox(0,0){$\beta=1$}}
\put(-7,3){\makebox(0,0){$\alpha=0$}}
\thicklines
\put(10,10){\line(1,0){50}}
\put(10,10){\line(0,1){50}}
\end{picture}
\end{center}
In this example, $(\alpha+\beta+\gamma)^2=2i$ and
$\alpha^2+\beta^2+\gamma^2=0$. Up to scale, therefore, (\ref{cube})  is the
correct equation.

The case of a regular dodecahedron is similar. Using the fact that a cube may
be inscribed in such a dodecahedron \cite{hilbertandcohn-vossen}, we may deduce
a particular projection:
\begin{center}
\begin{picture}(100,130)(-50,0)
\put(-450,-125){\epsfysize=1200pt \epsffile{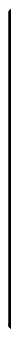}}
\put(-56,60){$\bullet$}
\put(-44,93){$\bullet$}
\put(-44,27){$\bullet$}
\put(-98,60){$\bullet$}
\put(-60,54){\scriptsize$0$}
\put(-94,54){\scriptsize$\beta=-1$}
\put(-37,26){\scriptsize$\displaystyle\gamma=\frac{\sqrt 5-1}{4}
                                            -\frac{\sqrt 5+1}{4}i$}
\put(-37,96){\scriptsize$\displaystyle\alpha=\frac{\sqrt 5-1}{4}
                                            +\frac{\sqrt 5+1}{4}i$}
\put(104,28){\scriptsize$1$}
\put(115,61){\tiny$\displaystyle\frac{\sqrt 5+1}{2}$}
\end{picture}
\end{center}
with $(\alpha+\beta+\gamma)^2=(7-3\sqrt 5)/2$ and
$\alpha^2+\beta^2+\gamma^2=(2-\sqrt 5)/2$. In this particular case,
$$(\alpha+\beta+\gamma)^2+(\sqrt 5-1)(\alpha^2+\beta^2+\gamma^2)=0.$$
Therefore, this is the correct equation in the general case. It may be used as
the basis of a ruler and compass construction of the general orthographic
projection of a regular dodecahedron.

It is interesting to note that if \underscore{all} the vertices of a Platonic
solid are orthographically projected to $z_1,z_2,\ldots,z_N\in{\Bbb C}$, then
\begin{equation}\label{platonic}(z_1+z_2+\cdots+z_N)^2=
  N(z_1{}^2+z_2{}^2+\cdots+z_N{}^2)\end{equation}
(compare (\ref{regularsimplex})). For a tetrahedron, this is just
equation~(\ref{tetrahedron}). To verify (\ref{platonic}) for the other Platonic
solids, firstly note that it is translation invariant. Therefore, it suffices
to impose $z_1+z_2+\cdots+z_N=0$ and show that
$z_1{}^2+z_2{}^2+\cdots+z_N{}^2=0$. The case of a cube now follows immediately
since its vertices may be grouped as two regular tetrahedra. The dodecahedral
case may be dealt with by grouping its vertices into five regular tetrahedra.
The regular octahedron is amenable to a similar trick but not the icosahedron.
Rather than resorting to direct computation, a uniform proof may be given as
follows. As before, assemble the vertices of the given solid $\Sigma$ as the
columns of a matrix~$A$, now of size~$3\times N$, and consider the moment
matrix~$M\equiv AA^t$. Observe that
$$\left\lgroup\begin{array}{ccc} 1&i&0\end{array}\right\rgroup M
  \left\lgroup\begin{array}c 1\\ i\\ 0\end{array}\right\rgroup
  =z_1{}^2+z_2{}^2+\cdots+z_N{}^2.$$
The moment matrix is positive definite and symmetric. In other words, it
defines a metric on ${\Bbb R}^3$, manifestly invariant under the symmetries
of~$\Sigma$. If $\Sigma$ is regular---or, more generally, enjoys the symmetries
of a regular solid (e.g.~a cuboctahedron or rhombicosidodecahedron)---then its
symmetry group acts irreducibly on~${\Bbb R}^3$. Thus, $M$ must be proportional
to the identity matrix and the result follows. For a general solid $\Sigma$,
the two complex numbers
$$\pm\sqrt{z_1{}^2+z_2{}^2+\cdots+z_N{}^2}$$
are the foci of the ellipse
$$\left\lgroup\begin{array}{ccc} x&y\end{array}\right\rgroup R
  \left\lgroup\begin{array}c x\\ y\end{array}\right\rgroup=1$$
where $R$ is the inverse of the quadratic form obtained by restricting M to the
plane of projection.

This reasoning also works in higher dimensions where it shows (as conjectured
to us by H.S.M. Coxeter) that the orthogonally projected images in the plane of
the $N$ vertices of any regular polytope, real or complex, will satisfy
equation (\ref{platonic}). Of course, this excludes regular polygons (whose
symmetry groups act reducibly except in dimension two) orthographic images
of which will satisfy (\ref{platonic}) if and only if the image is itself
regular. For polyhedra other than simplices, a quadratic equation such as
(\ref{platonic}) is no longer sufficient to characterise the orthogonal image
up to scale. In general, there will also be some linear relations. For a
non-degenerate $N$-gon there will be $N-n-1$ such relations. The simplest
example is a square in ${\Bbb R}^2$ which is characterised by the complex
equations
$$(\alpha+\beta+\gamma+\delta)^2=4(\alpha^2+\beta^2+\gamma^2+\delta^2)\quad
\mbox{ and }\quad\alpha+\gamma=\beta+\delta.$$

It is interesting to investigate further the relationship between a
non-degenerate simplex $\Sigma$ in ${\Bbb R}^n$ and its quadratic
form $Q=A^t(AA^t)^{-2}A$. Recall that $A$ is the $n\times(n+1)$ matrix whose
columns are the vertices of~$\Sigma$. There are several other formulae for or
characterisations of~$Q$. Let $S$ denote the $(n+1)\times(n+1)$ symmetric
matrix
$$\mbox{\huge $1$}\;
\raisebox{4pt}{$\displaystyle-\;\frac{1}{n+1}\left\lgroup\begin{array}{cccc}
                                     1&1&\cdots&1\\
                                     1&1&\cdots&1\\
                                     \vdots&\vdots&\ddots&\vdots\\
                                     1&1&\cdots&1\end{array}\right\rgroup$}.$$
It is the matrix of orthogonal projection in ${\Bbb R}^{n+1}$ in the direction 
of the vector~$e$. We maintain that $Q$ is characterised by the equations
$$QA^tA=S\quad\mbox{ and }\quad Qe=0.$$
Certainly, if these equations hold, then they are enough to determine $Q$
because the moment matrix $M\equiv A^tA$ has rank $n$ and $e$ is not in its
range. The second equation is evident and the first equation with $Q$ replaced
by $A^t(AA^t)^{-2}A$ and simplified reads
$$A^t(AA^t)^{-1}A=S.$$
To see that this holds it suffices to observe that it is clearly true after
post-multiplication by $A^t$ or~$e$. We may equally well characterise $Q$ by
means of the equations
$$A^tAQ=S\quad\mbox{ and }\quad Qe=0.$$
These equations relate $M$ and $Q$ geometrically---both matrices annihilate
$e$ whilst on the hyperplane orthogonal to $e$ they are mutually inverse. This
implies that $M$ and $Q$ are \underscore{generalised} \underscore{inverses}
\cite{penrose} of each other. Thus,
$$Q=M^\dagger=(A^tA)^\dagger=A^\dagger A^{\dagger t}$$
where $A^\dagger$ is the generalised inverse of~$A$. In this case,
$A^\dagger=A^t(AA^t)^{-1}$. This also shows how to compute $Q$
more directly in certain cases. The moment matrix $M$ has direct geometric
interpretation as the various inner products of the vectors
$a_1,a_2,\ldots,a_n,a_{n+1}$. In the case of a regular simplex, for example, we
know that $\|a_i\|^2$ is independent of $i$, that $\|a_i-a_j\|^2$ is
independent of $i\not= j$, and that $a_1+a_2+\ldots+a_n+a_{n+1}=0$. We may
deduce that, with a suitable overall scale, $M=S$. Since $S^\dagger=S$, it
follows that~$Q=S$. This is a direct derivation
of~(\ref{regularsimplex}).

It is clear geometrically that $M$ or, equivalently, $Q$ determines $\Sigma$
up to congruency. Alternatively, one can argue algebraically---it is easy to
check that if $A^tA=B^tB$, then $U=AA^t(BA^t)^{-1}$ is orthogonal and $A=UB$.
Therefore, the possible quadratic forms $Q$ which can arise give a natural
parametrisation of the non-degenerate simplices up to congruency. Choosing
a basepoint $\Sigma_0$ with corresponding matrix $A_0$, and mapping
$X\in{\mathrm{GL}}(n,{\Bbb R})$ to $X^{-1}A_0$ identifies the space of
non-degenerate simplices up to congruency with the homogeneous space
${\mathrm{GL}}(n,{\Bbb R})/{\mathrm O}(n)$. This homogeneous space may also be
identified with the space of positive definite $n\times n$ quadratic forms by
sending $X\in{\mathrm{GL}}(n,{\Bbb R})$ to~$XX^t$. The $(n+1)\times(n+1)$ 
quadratic form $Q$ corresponding to $X^{-1}A_0$ is given by
$A_0^\dagger XX^tA_0^{\dagger t}$. It follows that the general $Q$ which can
arise is characterised by the following two conditions:
\begin{itemize}
\item $Qe=0$ and only multiples of $e$ are in the kernel of~$Q$.
\item All other eigenvalues of $Q$ are positive.
\end{itemize}
It is also possible to repeat this analysis in pseudo-Euclidean spaces. The
only difference is that the condition that the non-zero eigenvalues of $Q$ be
positive is replaced by a condition on sign precisely reflecting the original
signature of the inner product.

Finally we should mention some possible applications. There is much current
interest in \underscore{computer} \underscore{vision}. In particular, there is
the problem of recognising a wire-frame object from its orthographic image. The
results we have described can be used as test on such an image, for example to
see whether a given image could be that of a cube or to keep track of a moving
shape. It is clear that such tests could be implemented quite efficiently.
Another possibility is in the manipulation of CADD\footnote{Computer Aided
Drafting and Design.} data. Rather than storing an image as an array of vectors
in ${\Bbb R}^3$, it may be sometimes be more efficient to store certain
tetrahedra within such an image by means of the corresponding quadratic form.
For orthographic imaging this may be preferable.

We would like to thank H.S.M.~Coxeter for drawing our attention to Hadwiger's
article, R.~Michaels and J.~Cofman for pointing out Gauss'
and Weisbach's work, and E.J.~Pitman for many useful conversations.

\small\renewcommand{\section}{\subsubsection}

\begin{tabbing}
Department of Pure Mathematics space \= \kill
Department of Pure Mathematics       \> Mathematical Institute\\
University of Adelaide               \> 24-29 Saint Giles'\\
South AUSTRALIA 5005                 \> Oxford OX1 3LB\\
                                     \> ENGLAND
\end{tabbing}
\end{document}